\documentclass[11pt,thmsa,arabtex]{article}
\usepackage{amssymb}
\usepackage{amsfonts}
\usepackage{amsmath}
\usepackage{graphicx}
\usepackage{hyperref}
\usepackage{cite}
\usepackage{ntheorem}
\usepackage{epstopdf}

\newtheorem*{theorem-non}{Theorem}

\newtheorem{case}{Case}[section]

\newtheorem{definition}{Definition}[section]
\newtheorem{notation}{Notation}[section]

\newtheorem{example}{Example}[section]

\numberwithin{theorem}{section}
\numberwithin{equation}{section}

\begin{document}

\title{{\textbf{Spacelike and Timelike Admissible Smarandache Curves in Pseudo-Galilean Space}}}
\author{\textbf{M. Khalifa Saad}\thanks{
~E-mail address:~mohamed\_khalifa77@science.sohag.edu.eg} \\
{\small Mathematics Dept., Faculty of Science, Sohag University, 82524 Sohag, Egypt}}
\date{}
\maketitle

\textbf{Abstract.} In this paper, space and timelike admissible Smarandache curves in the pseudo-Galilean space $G_{3}^{1}$ are investigated. Also, Smarandache curves of the position vector of space and timelike arbitrary curve and some of its special curves in $G_{3}^{1}$ are obtained. To confirm our main results, some examples are given and illustrated.

\bigskip \textbf{\emph{M.S.C. 2010}:} 53A35, 53B30.


\section{Introduction}
In recent years, researchers have begun to investigate curves and surfaces in the Galilean space and thereafter pseudo-Galilean space $G_{3}$ and $G_{3}^{1}$. In the study of the fundamental theory and the characterizations of space curves, the corresponding relations between the curves are the very interesting and important problem.\\
It is known that a Smarandache geometry is a geometry which has at least one Smarandache denied axiom \cite{1}. An axiom is said to be Smarandache denied, if it behaves in at least two different ways within the same space. Smarandache geometries are connected with the theory of relativity and the parallel universes. Smarandasche curves are the objects of Smarandache geometry. By definition, if the position vector of a curve $\delta$ is composed by Frenet frame's vectors of another curve $\beta$, then the curve $\delta$ is called a Smarandache curve \cite{2}. Smarandache curves have been investigated by some differential geometers (see for example, \cite{2,3}). M. Turgut and S. Yilmaz  defined a special case of such curves and call it Smarandache $TB_{2}$ curves in the space $E_{1}^{4}$ \cite{2}. They studied special Smarandache curves which are defined by the tangent and second binormal vector fields. In \cite{3}, the author introduced some special Smarandache curves in the Euclidean space. He studied Frenet-Serret invariants of a special case.\\
In the field of computer aided design and computer graphics, helices can be used for the tool path description, the simulation of kinematic motion or the design of highways, etc. \cite{4}. The main feature of general helix is that the tangent makes a constant angle with a fixed straight line which is called the axis of the general helix. A classical result stated by Lancret in 1802 and first proved by de Saint Venant in 1845 says that: A necessary and sufficient condition that a curve be a general helix is that the ratio ($\kappa / \tau$) is constant along the curve, where $\kappa$ and $\tau$ denote the curvature and the torsion, respectively. Also, the helix is also known as circular helix or W-curve which is a special case of the general helix \cite{5}.\\
Salkowski (resp. Anti-Salkowski) curves in Euclidean space are generally known as family of curves with constant curvature (resp. torsion) but nonconstant torsion (resp. curvature) with an explicit parametrization. They were defined in an earlier paper \cite{6}.\\
In this paper, we obtain Smarandache curves for a position vector of an arbitrary curve in $G_{3}^{1}$ and some of its special curves (helix, circular helix, Salkowski and Anti-Salkowski curves). In other words, according to Frenet frame {$\textbf{e}_{1}$, $\textbf{e}_{2}$, $\textbf{e}_{3}$} of the considered curves in the pseudo-Galilean space $G_{3}^{1}$, the meant Smarandache curves $\textbf{e}_{1}\textbf{e}_{2}$, $\textbf{e}_{1} \textbf{e}_{3}$ and $\textbf{e}_{1} \textbf{e}_{2} \textbf{e}_{3}$ are obtained.
To the best of author's knowledge, Smarandache curves have not been presented in the pseudo-Galilean geometry in depth. Thus, the study is proposed to serve such a need.

\section{Basic notions and properties}

In this section, let us first recall basic notions from pseudo-Galilean
geometry \cite{7,8,9,10,11}. In the inhomogeneous affine coordinates for points and vectors
(point pairs) the similarity group $H_{8}$ of $G_{3}^{1}$ has the following
form
\begin{align}
\bar{x}& =a+b.x,  \notag \\
\bar{y}& =c+d.x+r.\cosh \theta .y+r.\sinh \theta .z,  \notag \\
\bar{z}& =e+f.x+r.\sinh \theta .y+r.\cosh \theta .z,
\end{align}%
where $a,b,c,d,e,f,r$ and $\theta $ are real numbers.
Particularly, for $b=r=1,$ the group $(2.1)$ becomes the group $B_{6}\subset
H_{8}$ of isometries (proper motions) of the pseudo-Galilean space $G_{3}^{1}
$. The motion group leaves invariant the absolute figure and defines the
other invariants of this geometry. It has the following form
\begin{align}
\bar{x}& =a+x,  \notag \\
\bar{y}& =c+d.x+\cosh \theta .y+\sinh \theta .z,  \notag \\
\bar{z}& =e+f.x+\sinh \theta .y+\cosh \theta .z.
\end{align}%
According to the motion group in the pseudo-Galilean space, there are
non-isotropic vectors $A(A_{1},A_{2},A_{3})$ (for which holds $A_{1}\neq 0$)
and four types of isotropic vectors: spacelike ($A_{1}=0,$ $%
A_{2}^{2}-A_{3}^{2}>0$), timelike ($A_{1}=0,$ $A_{2}^{2}-A_{3}^{2}<0$) and
two types of lightlike vectors ($A_{1}=0,A_{2}=\pm A_{3}$). The scalar
product of two vectors $u=(u_{1},u_{2},u_{3})$ and $v=(v_{1},v_{2},v_{3})$
in $G_{3}^{1}$ is defined by
\begin{equation*}
\left\langle u,v\right\rangle _{G_{3}^{1}}=\left\{
\begin{array}{c}
u_{1}v_{1},\text{ \ \ \ \ \ \ \ \ \ \ \ \ \ \ \ \ \ \ \ \ \ \ \ \ if }%
u_{1}\neq 0\text{ or }v_{1}\neq 0, \\
u_{2}v_{2}-u_{3}v_{3}\text{ \ \ \ \ \ \ \ \ \ \ \ \ \ \ \ \ \ \ \ \ if \ }%
u_{1}=0\text{ and }v_{1}=0.%
\end{array}
\right\}
\end{equation*}
We introduce a pseudo-Galilean cross product in the following way
\begin{equation*}
u\times _{G_{3}^{1}}v=\left\vert
\begin{array}{ccc}
0 & -j & k \\
u_{1} & u_{2} & u_{3} \\
v_{1} & v_{2} & v_{3}%
\end{array}%
\right\vert ,
\end{equation*}%
where $j=(0,1,0)$ and $k=(0,0,1)$ are unit spacelike and timelike vectors,
respectively. Let us recall basic facts about curves in $G_{3}^{1}$, that were introduced in \cite{7,8,9}.

A curve $\gamma (s)=(x(s),y(s),z(s))$ is called an admissible curve if it has no
inflection points $(\dot{\gamma}\times \ddot{\gamma}\neq 0)$ and no
isotropic tangents $(\dot{x}\neq 0)$ or normals whose projections on the
absolute plane would be lightlike vectors $(\dot{y}\neq \pm \dot{z})$. An
admissible curve in $G_{3}^{1}$ is an analogue of a regular curve in
Euclidean space \cite{8}.

For an admissible curve $\gamma :I\subseteq \mathbb{R}\rightarrow G_{3}^{1},$
the curvature $\kappa (s)$ and torsion $\tau (s)$ are defined by%
\begin{equation}
\kappa (s)=\frac{\sqrt{\left\vert \ddot{y}(s)^{2}-\ddot{z}(s)^{2}\right\vert
}}{(\dot{x}(s))^{2}},\text{ }\tau (s)=\frac{\ddot{y}(s)\dddot{z}(s)-\dddot{y}%
(s)\ddot{z}(s)}{\left\vert \dot{x}(s)\right\vert ^{5}\cdot \kappa ^{2}(s)},%
\text{\ }
\end{equation}%
expressed in components. Hence, for an admissible curve $\gamma :I\subseteq
\mathbb{R}\rightarrow G_{3}^{1}$ parameterized by the arc length $s$ with
differential form $ds=dx$, given by
\begin{equation}
\gamma (x)=(x,y(x),z(x)),
\end{equation}%
the formulas $(2.3)$ have the following form
\begin{equation}
\kappa (x)=\sqrt{\left\vert y^{^{\prime \prime }}(x)^{2}-z^{^{\prime \prime
}}(x)^{2}\right\vert },\text{ }\tau (x)=\frac{y^{^{\prime \prime
}}(x)z^{^{\prime \prime \prime }}(x)-y^{^{\prime \prime \prime
}}(x)z^{^{\prime \prime }}(x)}{\kappa ^{2}(x)}.
\end{equation}%
The associated trihedron is given by
\begin{align}
\mathbf{e}_{1}& =\gamma ^{\prime }(x)=(1,y^{^{\prime }}(x),z^{^{\prime
}}(x)),  \notag \\
\mathbf{e}_{2}& =\frac{1}{\kappa (x)}\gamma ^{^{\prime \prime }}(x)=\frac{1}{%
\kappa (x)}(0,y^{^{\prime \prime }}(x),z^{^{\prime \prime }}(x)),  \notag \\
\mathbf{e}_{3}& =\frac{1}{\kappa (x)}(0,\epsilon z^{^{\prime \prime
}}(x),\epsilon y^{^{\prime \prime }}(x)),
\end{align}%
where $\epsilon =+1$ or $\epsilon =-1$, chosen by criterion Det$(e_{1},e_{2},e_{3})=1$, that means
\begin{equation*}
\left\vert y^{^{\prime \prime }}(x)^{2}-z^{^{\prime \prime
}}(x)^{2}\right\vert =\epsilon (y^{^{\prime \prime }}(x)^{2}-z^{^{\prime
\prime }}(x)^{2})\text{.}
\end{equation*}
The curve $\gamma $ given by $(2.4)$ is timelike (resp. spacelike) if $%
\mathbf{e}_{2}(s)$ is a spacelike (resp. timelike) vector. The principal
normal vector or simply normal is spacelike if $\epsilon =+1$ and timelike
if $\epsilon =-1$. For derivatives of the tangent $\mathbf{e}_{1}$, normal $%
\mathbf{e}_{2}$ and binormal $\mathbf{e}_{3}$ vector fields, the following
Frenet formulas in $G_{3}^{1}$ hold:
\begin{align}
\mathbf{e}_{1}^{\prime }(x)& =\kappa (x)\mathbf{e}_{2}(x),  \notag \\
\mathbf{e}_{2}^{\prime }(x)& =\tau (x)\mathbf{e}_{3}(x),  \notag \\
\mathbf{e}_{3}^{\prime }(x)& =\tau (x)\mathbf{e}_{2}(x).
\end{align}
From $(2.5)$ and $(2.6)$, we have the following important relation that is true in Galilean and pseudo-Galilean spaces \cite{11,12,13}
\begin{equation*}
\eta ^{\prime \prime \prime }(s)=\kappa ^{\prime }(s)\mathbf{N}(s)+\kappa
(s)\tau (s)\mathbf{B}(s).
\end{equation*}
In \cite{2} authors introduced:
\begin{definition}
A regular curve in Minkowski space-time, whose position vector is composed by Frenet frame vectors on another regular curve, is called a Smarandache curve.
\end{definition}
In the light of the above definition, we adapt it to admissible curves in the pseudo-Galilean space as
follows:
\begin{definition}
let $\eta =\eta (s)$ be an admissible curve in $G_{3}^{1}$ and $\{\mathbf{e}_{1},
\mathbf{e}_{2},\mathbf{e}_{3}\}$ be its moving Frenet frame. Smarandache $\mathbf{e}_{1}\mathbf{e}_{2}
,\mathbf{e}_{1}\mathbf{e}_{3}$ and $\mathbf{e}_{1}\mathbf{e}_{2}\mathbf{e}_{3}$ curves are respectively, defined by
\begin{eqnarray}
\eta _{\mathbf{e}_{1}\mathbf{e}_{2}} &=&\frac{\mathbf{e}_{1}+\mathbf{e}_{2}}{\left\Vert \mathbf{e}_{1}+%
\mathbf{e}_{2}\right\Vert },\ \ \   \notag \\
\eta _{\mathbf{e}_{1}\mathbf{e}_{3}} &=&\frac{\mathbf{e}_{1}+\mathbf{e}_{3}}{\left\Vert \mathbf{e}_{1}+%
\mathbf{e}_{3}\right\Vert },  \notag \\
\eta _{\mathbf{e}_{1}\mathbf{e}_{2}\mathbf{e}_{3}} &=&\frac{\mathbf{e}_{1}+\mathbf{\mathbf{e}_{2}+\mathbf{e}_{3}}}{\left\Vert \mathbf{\mathbf{e}_{1}
}+\mathbf{\mathbf{e}_{2}+\mathbf{e}_{3}}\right\Vert }.
\end{eqnarray}
\end{definition}

\section{Smarandache curves of an arbitrary curve in $G_{3}^{1}$}

In the light of which introduced in the Galilean 3-space $G_{3}$ by \cite{3}, we introduce the position vectors of \emph{spacelike} and \emph{timelike} arbitrary curves with curvature $\kappa (s)$ and torsion $\tau (s)$ in the pseudo-Galilean space $G_{3}^{1}$ and then calculate their Smarandache curves.\\
Let us start with an arbitrary curve $\textbf{r}(s)$ in $G_{3}^{1}$, so we get\\
\begin{case}
$\textbf{r}(s)$ is \emph{spacelike}:
\begin{equation}
\mathbf{r}(s)=\left( s,-\int \left( \int \sinh \left( \int \tau (s)\,ds\right) \kappa
(s)\,ds\right) \,ds,\int \left( \int \cosh \left( \int \tau (s)\,ds\right)
\kappa (s)\,ds\right) \,ds\right).
\end{equation}
The derivatives of this curve are respectively, given by
\begin{equation*}
\mathbf{r}^{\prime }(s)=\left( 1,-\int \sinh \left( \int \tau (s)\,ds\right) \kappa (s)\,\,ds,\int
\cosh \left( \int \tau (s)\,ds\right) \kappa (s)\,ds\right),
\end{equation*}
\begin{equation*}
\mathbf{r}^{\prime \prime }(s)=\left( 0,-\sinh \left( \int \tau (s)\,ds\right) \kappa (s),\cosh \left( \int
\tau (s)\,ds\right) \kappa (s)\right),
\end{equation*}
\begin{equation}
\mathbf{r}^{\prime \prime \prime }(s)=\left(
\begin{array}{c}
0,-\kappa ^{\prime }\sinh \left( \int \tau (s)\,ds\right) -\cosh \left( \int
\tau (s)\,ds\right) \kappa (s)\tau (s), \\
\cosh \left( \int \tau (s)\,ds\right) \kappa ^{\prime }+\sinh \left( \int
\tau (s)\,ds\right) \kappa (s)\tau (s)
\end{array}
\right).
\end{equation}
The frame vector fields of $\textbf{r}$ are as follows
\begin{equation*}
(\mathbf{e}_{1})_{\textbf{r}}=\left( 1,-\int \sinh \left( \int \tau (s)\,ds\right) \kappa (s)ds,\int
\cosh \left( \int \tau (s)\,ds\right) \kappa (s)\,ds\right),
\end{equation*}
\begin{equation*}
(\mathbf{e}_{2})_{\textbf{r}}=\left( 0,-\sinh \left( \int \tau (s)\,ds\right) ,\cosh \left( \int \tau
(s)\,ds\right) \right),
\end{equation*}
\begin{equation}
(\mathbf{e}_{3})_{\textbf{r}}=\left( 0,-\cosh \left( \int \tau (s)\,ds\right) ,\sinh \left( \int \tau
(s)\,ds\right) \right).
\end{equation}
By Definition (2.2), the $\textbf{e}_{1}\textbf{e}_{2}$, $\textbf{e}_{1}\textbf{e}_{3}$ and $\textbf{e}_{1}\textbf{e}_{2}\textbf{e}_{3}$ Smarandache curves of $r$ are respectively, written as
\begin{equation*}
\mathbf{r}_{\mathbf{e}_{1}\mathbf{e}_{2}}=\left(
\begin{array}{c}
1,-\int \sinh \left( \int \tau (s)\,ds\right) \kappa (s)\,ds-\sinh \left(
\int \tau (s)\,ds\right) , \\
\cosh \left( \int \tau (s)\,ds\right) +\int \cosh \left( \int \tau
(s)\,ds\right) \kappa (s)\,ds%
\end{array}%
\right),
\end{equation*}
\begin{equation*}
\mathbf{r}_{\mathbf{e}_{1}\mathbf{e}_{3}}=\left(
\begin{array}{c}
1,-\cosh \left( \int \tau (s)\,ds\right) -\int \sinh \left( \int \tau
(s)\,ds\right) \kappa (s)\,ds, \\
\int \cosh \left( \int \tau (s)\,ds\right) \kappa (s)ds+\sinh \left( \int
\tau (s)\,ds\right)
\end{array}%
\right),
\end{equation*}
\begin{equation}
\mathbf{r}_{\mathbf{e}_{1}\mathbf{e}_{2}\mathbf{e}_{3}}=\left(
\begin{array}{c}
1,-e^{\int \tau \lbrack s]\,ds}-\int \sinh \left( \int \tau (s)\,ds\right)
\kappa (s)\,ds, \\
e^{\int \tau \lbrack s]\,ds}+\int \cosh \left( \int \tau (s)\,ds\right)
\kappa (s)ds%
\end{array}%
\right).
\end{equation}
\end{case}
\begin{case}
$\textbf{r}(s)$ is \emph{timelike}:
\begin{equation}
\mathbf{r}(s)=\left( s,\int \left( \int \cosh \left( \int \tau (s)\,ds\right) \kappa
(s)\,ds\right) \,ds,\int \left( \int \sinh \left( \int \tau (s)\,ds\right)
\kappa (s)\,ds\right) \,ds\right).
\end{equation}
So, the derivatives of $\textbf{r}(s)$ are
\begin{equation*}
\mathbf{r}^{\prime }(s)=\left( 1,\int \cosh \left( \int \tau (s)\,ds\right) \kappa (s)\,ds,\int
\sinh \left( \int \tau (s)\,ds\right) \kappa (s)\,ds\right),
\end{equation*}
\begin{equation*}
\mathbf{r}^{\prime \prime }(s)=\left( 0,\cosh \left( \int \tau (s)\,ds\right) \kappa (s)\,,\sinh \left(
\int \tau (s)\,ds\right) \kappa (s)\right),
\end{equation*}
\begin{equation}
\mathbf{r}^{\prime \prime \prime }(s)=\left(
\begin{array}{c}
0,\cosh \left( \int \tau (s)\,ds\right) \kappa ^{\prime }+\sinh \left( \int
\tau (s)\,ds\right) \kappa (s)\tau (s), \\
\kappa ^{\prime }\sinh \left( \int \tau (s)\,ds\right) +\cosh \left( \int
\tau (s)\,ds\right) \kappa (s)\tau (s)%
\end{array}%
\right).
\end{equation}
And the frame vector fields are as follows
\begin{equation*}
(\mathbf{e}_{1})_{\textbf{r}}=\left( 1,\int \cosh \left( \int \tau (s)\,ds\right) \kappa (s)\,ds,\int
\sinh \left( \int \tau (s)\,ds\right) \kappa (s)\,ds\right),
\end{equation*}
\begin{equation*}
(\mathbf{e}_{2})_{\textbf{r}}=\left( 0,\cosh \left( \int \tau (s)\,ds\right) ,\sinh \left( \int \tau
(s)\,ds\right) \right),
\end{equation*}
\begin{equation}
(\mathbf{e}_{3})_{\textbf{r}}=\left( 0,\sinh \left( \int \tau (s)\,ds\right) ,\cosh \left( \int \tau
(s)\,ds\right) \right).
\end{equation}
Hence, the Smarandache curves are
\begin{equation*}
\mathbf{r}_{\mathbf{e}_{1}\mathbf{e}_{2}}=\left(
\begin{array}{c}
1,\cosh \left( \int \tau (s)\,ds\right) +\int \cosh \left( \int \tau
(s)\,ds\right) \kappa (s)\,ds, \\
\int \sinh \left( \int \tau (s)\,ds\right) \kappa (s)\,ds+\sinh \left( \int
\tau (s)\,ds\right)
\end{array}%
\right),
\end{equation*}
\begin{equation*}
\mathbf{r}_{\mathbf{e}_{1}\mathbf{e}_{3}}=\left(
\begin{array}{c}
1,\int \cosh \left( \int \tau (s)\,ds\right) \kappa (s)\,ds+\sinh \left(
\int \tau (s)\,ds\right) , \\
\cosh \left( \int \tau (s)\,ds\right) +\int \sinh \left( \int \tau
(s)\,ds\right) \kappa (s)\,ds%
\end{array}%
\right),
\end{equation*}
\begin{equation}
\mathbf{r}_{\mathbf{e}_{1}\mathbf{e}_{2}\mathbf{e}_{3}}=\left(
\begin{array}{c}
1,e^{\int \tau (s)ds}+\int \cosh \left( \int \tau (s)\,ds\right) \kappa
(s)ds, \\
e^{\int \tau (s)\,ds}+\int \sinh \left( \int \tau (s)\,ds\right) \kappa (s)ds%
\end{array}%
\right).
\end{equation}
\end{case}
\section{Smarandache curves of some special curves in $G_{3}^{1}$}
\subsection{Smarandache curves of a general helix}

Let $\alpha(s)$ be an admissible general helix in $G_{3}^{1}$ with ($\tau /\kappa=m=const.$), we have\\
\textbf{Case 4.1.1} \ $\alpha(s)$ is \emph{spacelike}:
\begin{equation}
\alpha (s)=\left( s,-\frac{1}{m}\int \cosh \left( m\int \kappa (s)\,ds\right) \,ds,%
\frac{1}{m}\int \sinh \left( m\int \kappa (s)\,ds\right) \,ds\right).
\end{equation}
Then $\alpha ^{\prime }, \alpha ^{\prime \prime }, \alpha ^{\prime \prime \prime }$ for this curve are respectively, expressed as
\begin{equation*}
\alpha ^{\prime }(s)=\left( 1,-\frac{1}{m}\cosh \left( m\int \kappa (s)\,ds\right) ,\frac{1}{m}%
\sinh \left( m\int \kappa (s)\,ds\right) \right),
\end{equation*}
\begin{equation*}
\alpha ^{\prime \prime }(s)=\left( 0,-\sinh \left( m\int \kappa (s)\,ds\right) \kappa (s),\cosh \left(
m\int \kappa (s)\,ds\right) \kappa (s)\right),
\end{equation*}
\begin{equation}
\alpha ^{\prime \prime \prime }(s)=\left(
\begin{array}{c}
0,-\kappa ^{\prime }\sinh \left( m\int \kappa (s)\,ds\right) -m\cosh \left(
m\int \kappa (s)\,ds\right) \kappa ^{2}(s), \\
\cosh \left( m\int \kappa (s)\,ds\right) \kappa ^{\prime }+m\sinh \left(
m\int \kappa (s)\,ds\right) \kappa ^{2}(s)%
\end{array}%
\right).
\end{equation}
The moving Frenet vectors of $\alpha (s)$ are given by
\begin{equation*}
(\mathbf{e}_{1})_{\alpha}=\left( 1,-\frac{1}{m}\cosh \left( m\int \kappa (s)\,ds\right) ,\frac{1}{m}%
\sinh \left( m\int \kappa (s)\,ds\right) \right),
\end{equation*}
\begin{equation*}
(\mathbf{e}_{2})_{\alpha}=\left( 0,-\sinh \left( m\int \kappa (s)\,ds\right) ,\cosh \left( m\int
\kappa (s)\,ds\right) \right),
\end{equation*}
\begin{equation}
(\mathbf{e}_{3})_{\alpha}=\left( 0,-\cosh \left( m\int \kappa (s)\,ds\right) ,\sinh \left( m\int
\kappa (s)\,ds\right) \right).
\end{equation}
From which, Smarandache curves are given by
\begin{equation*}
\alpha_{\mathbf{e}_{1}\mathbf{e}_{2}}=\left(
\begin{array}{c}
1,-\frac{1}{m}\cosh \left( m\int \kappa (s)\,ds\right) +m\sinh \left( m\int
\kappa (s)\,ds\right) , \\
\cosh \left( m\int \kappa (s)\,ds\right) +\frac{1}{m}\sinh \left( m\int
\kappa (s)\,ds\right)
\end{array}%
\right),
\end{equation*}
\begin{equation*}
\alpha_{\mathbf{e}_{1}\mathbf{e}_{3}}=\left( 1,-\frac{1}{m}(1+m)\cosh \left( m\int \kappa (s)\,ds\right) ,\frac{1}{%
m}(1+m)\sinh \left( m\int \kappa (s)\,ds\right) \right),
\end{equation*}
\begin{equation}
\alpha_{\mathbf{e}_{1}\mathbf{e}_{2}\mathbf{e}_{3}}=\left(
\begin{array}{c}
1,-\frac{1}{m}(1+m)\cosh \left( m\int \kappa (s)\,ds\right) +m\sinh \left(
m\int \kappa (s)\,ds\right) , \\
e^{m\int \kappa (s)\,ds}+\frac{1}{m}\sinh \left( m\int \kappa (s)\,ds\right)
\end{array}
\right).
\end{equation}\\
\textbf{Case 4.1.2} \ $\alpha(s)$ is \emph{timelike}:
\begin{equation}
\alpha (s)=\left( s,\frac{1}{m}\int \sinh \left( m\int \kappa (s)\,ds\right) \,ds,\frac{%
1}{m}\int \cosh \left( m\int \kappa (s)\,ds\right) \,ds\right).
\end{equation}
So, $\alpha ^{\prime }, \alpha ^{\prime \prime }, \alpha ^{\prime \prime \prime }$ are respectively,
\begin{equation*}
\alpha ^{\prime }(s)=\left( 1,\frac{1}{m}\sinh \left( m\int \kappa (s)\,ds\right) ,\frac{1}{m}%
\cosh \left( m\int \kappa (s)\,ds\right) \right),
\end{equation*}
\begin{equation*}
\alpha ^{\prime \prime }(s)=\left( 0,\cosh \left( m\int \kappa (s)\,ds\right) \kappa (s),\sinh \left(
m\int \kappa (s)\,ds\right) \kappa (s)\right),
\end{equation*}
\begin{equation}
\alpha ^{\prime \prime \prime }(s)=\left(
\begin{array}{c}
0,\kappa ^{\prime }\cosh \left( m\int \kappa (s)\,ds\right) +m\sinh \left(
m\int \kappa (s)\,ds\right) \kappa ^{2}(s), \\
\sinh \left( m\int \kappa (s)\,ds\right) \kappa ^{\prime }+m\cosh \left(
m\int \kappa (s)\,ds\right) \kappa ^{2}(s)%
\end{array}%
\right).
\end{equation}
The Frenet vectors of $\alpha (s)$ are given by
\begin{equation*}
(\mathbf{e}_{1})_{\alpha}=\left( 1,\frac{1}{m}\sinh \left( m\int \kappa (s)\,ds\right) ,\frac{1}{m}%
\cosh \left( m\int \kappa (s)\,ds\right) \right),
\end{equation*}
\begin{equation*}
(\mathbf{e}_{2})_{\alpha}=\left( 0,\cosh \left( m\int \kappa (s)\,ds\right) ,\sinh \left( m\int \kappa
(s)\,ds\right) \right),
\end{equation*}
\begin{equation}
(\mathbf{e}_{3})_{\alpha}=\left( 0,\sinh \left( m\int \kappa (s)\,ds\right) ,\cosh \left( m\int \kappa
(s)\,ds\right) \right).
\end{equation}
So, Smarandache curves of $\alpha$ are as follows
\begin{equation*}
\alpha_{\mathbf{e}_{1}\mathbf{e}_{2}}=\left(
\begin{array}{c}
1,\frac{1}{m}\sinh \left( m\int \kappa (s)\,ds\right) +\cosh \left( m\int
\kappa (s)\,ds\right) , \\
\frac{1}{m}\cosh \left( m\int \kappa (s)\,ds\right) +\sinh \left( m\int
\kappa (s)\,ds\right)
\end{array}%
\right),
\end{equation*}
\begin{equation*}
\alpha_{\mathbf{e}_{1}\mathbf{e}_{3}}=\left( 1,\frac{1}{m}(1+m)\sinh \left( m\int \kappa (s)\,ds\right) ,\frac{1}{m%
}(1+m)\cosh \left( m\int \kappa (s)\,ds\right) \right),
\end{equation*}
\begin{equation}
\alpha_{\mathbf{e}_{1}\mathbf{e}_{2}\mathbf{e}_{3}}=\left(
\begin{array}{c}
1,e^{m\int \kappa (s)\,ds}+\frac{1}{m}\sinh \left( m\int \kappa
(s)\,ds\right) , \\
\frac{1}{m}(1+m)\cosh \left( m\int \kappa (s)\,ds\right) +\sinh \left( m\int
\kappa (s)\,ds\right) ,%
\end{array}%
\right).
\end{equation}

\subsection{Smarandache curves of a circular helix}

Let $\beta(s)$ be an admissible circular helix in $G_{3}^{1}$ with ($\tau=a=const., \kappa=b=const.$), we have\\
\textbf{Case 4.2.1} \ $\beta(s)$ is \emph{spacelike}:
\begin{equation}
\beta (s)=\left( s,a\int \left( \int \sinh (bs)\,ds\right) \,ds,a\int \left( \int
\cosh (bs)\,ds\right) \,ds\right).
\end{equation}
For this curve, we have
\begin{equation*}
\beta ^{\prime }(s)=\left( 1,\frac{a}{b}\cosh (bs),\frac{a}{b}\sinh (bs)\right),
\end{equation*}
\begin{equation*}
\beta ^{\prime \prime }(s)=\left( 1,\frac{a}{b}\cosh (bs),\frac{a}{b}\sinh (bs)\right),
\end{equation*}
\begin{equation}
\beta ^{\prime \prime \prime }(s)=\left( 0,a\sinh (bs),a\cosh (bs)\right).
\end{equation}
Making necessary calculations from above, we have
\begin{equation*}
(\mathbf{e}_{1})_{\beta}=\left( 1,\frac{a}{b}\cosh (bs),\frac{a}{b}\sinh (bs)\right),
\end{equation*}
\begin{equation*}
(\mathbf{e}_{2})_{\beta}=\left( 0,\sinh (bs),\cosh (bs)\right),
\end{equation*}
\begin{equation}
(\mathbf{e}_{3})_{\beta}=\left( 0,-\cosh (bs),-\sinh (bs)\right).
\end{equation}
Considering the last Frenet vectors, the $\textbf{e}_{1}\textbf{e}_{2}$, $\textbf{e}_{1}\textbf{e}_{3}$ and $\textbf{e}_{1}\textbf{e}_{2}\textbf{e}_{3}$ Smarandache curves of $\beta$ are respectively, as follows
\begin{equation*}
\beta_{\mathbf{e}_{1}\mathbf{e}_{2}}=\left( 1,\frac{a}{b}\cosh (bs)+\sinh (bs),\cosh (bs)+\frac{a}{b}\sinh
(bs)\right),
\end{equation*}
\begin{equation*}
\beta_{\mathbf{e}_{1}\mathbf{e}_{3}}=\left( 1,\frac{(a-b)}{b}\cosh (bs),\frac{(a-b)}{b}\sinh (bs)\right),
\end{equation*}
\begin{equation}
\beta_{\mathbf{e}_{1}\mathbf{e}_{2}\mathbf{e}_{3}}=\left( 1,\left( \frac{a}{b}-1\right) \cosh (bs)+\sinh (bs),\cosh (bs)+\frac{(a-b)}{b}\sinh (bs)\right).
\end{equation}
\textbf{Case 4.2.2} \ $\beta(s)$ is \emph{timelike}:
\begin{equation}
\beta (s)=\left( s,-a\int \left( \int \cosh (bs)\,ds\right) \,ds,a\int \left( \int
\sinh (bs)\,ds\right) \,ds\right).
\end{equation}
For $\beta(s)$, we have
\begin{equation*}
\beta ^{\prime }(s)=\left( 1,-\frac{a}{b}\sinh (bs),\frac{a}{b}\cosh (bs)\right),
\end{equation*}
\begin{equation*}
\beta ^{\prime \prime }(s)=\left( 0,-a\cosh (bs),a\sinh (bs)\right),
\end{equation*}
\begin{equation}
\beta ^{\prime \prime \prime }(s)=\left( 0,-ab\sinh (bs),ab\cosh (bs)\right).
\end{equation}
The Frenet frame of $\beta$ is
\begin{equation*}
(\mathbf{e}_{1})_{\beta}=\left( 1,-\frac{a}{b}\sinh (bs),\frac{a}{b}\cosh (bs)\right),
\end{equation*}
\begin{equation*}
(\mathbf{e}_{2})_{\beta}=\left( 0,-\cosh (bs),\sinh (bs)\right),
\end{equation*}
\begin{equation}
(\mathbf{e}_{3})_{\beta}=\left( 0,\sinh (bs),-\cosh (bs)\right).
\end{equation}
Thus the Smarandache curves of $\beta$ are respectively, given by
\begin{equation*}
\beta_{\mathbf{e}_{1}\mathbf{e}_{2}}=\left( 1,-\frac{1}{b}\left( b\cosh (bs)+a\sinh (bs)\right) ,\frac{a}{b}\cosh
(bs)+\sinh (bs)\right),
\end{equation*}
\begin{equation*}
\beta_{\mathbf{e}_{1}\mathbf{e}_{3}}=\left( 1,-\frac{(a+b)}{b}\sinh (bs),\frac{(a+b)}{b}\cosh (bs)\right),
\end{equation*}
\begin{equation}
\beta_{\mathbf{e}_{1}\mathbf{e}_{2}\mathbf{e}_{3}}=\left( 1,-\frac{1}{b}\left( be^{bs}+a\sinh (bs)\right) ,\frac{(a+b)}{b}\cosh
(bs)+\sinh (bs)\right).
\end{equation}

\subsection{Smarandache curves of Salkowski curve}

Let $\gamma (s)$ be a Salkowski curve in $G_{3}^{1}$ with ($\tau=\tau (s), \kappa=a=const.$)\\
\textbf{Case 4.3.1} \ $\gamma(s)$ is \emph{spacelike}:
\begin{equation}
\gamma (s)=\left( s,-a\int \left( \int \sinh \left( \int \tau (s)\,ds\right)
\,ds\,\right) ds,a\int \left( \int \cosh \left( \int \tau (s)\,ds\right)
\,ds\,\right) ds\right).
\end{equation}
If we differentiate this equation three times, one can obtain
\begin{equation*}
\gamma ^{\prime }(s)=\left( 1,-a\int \sinh \left( \int \tau (s)\,ds\right) \,ds,a\int \cosh
\left( \int \tau (s)\,ds\right) \,ds\right),
\end{equation*}
\begin{equation*}
\gamma ^{\prime \prime }(s)=\left( 0,-a\sinh \left( \int \tau (s)\,ds\right) ,a\cosh \left( \int \tau
(s)\,ds\right) \right),
\end{equation*}
\begin{equation}
\gamma ^{\prime \prime \prime }(s)=\left( 0,-a\cosh \left( \int \tau (s)\,ds\right) \tau (s),a\sinh \left( \int
\tau (s)\,ds\right) \tau (s)\right).
\end{equation}
In addition to that, the tangent, principal normal and binormal vectors of $\gamma$ are in the following forms
\begin{equation*}
(\mathbf{e}_{1})_{\gamma}=\left( 1,-a\int \sinh \left( \int \tau (s)\,ds\right) \,ds,a\int \cosh
\left( \int \tau (s)\,ds\right) \,ds\right),
\end{equation*}
\begin{equation*}
(\mathbf{e}_{2})_{\gamma}=\left( 0,-\sinh \left( \int \tau (s)\,ds\right) ,\cosh \left( \int \tau
(s)\,ds\right) \right),
\end{equation*}
\begin{equation}
(\mathbf{e}_{3})_{\gamma}=\left( 0,-\cosh \left( \int \tau (s)\,ds\right) ,\sinh \left( \int \tau
(s)\,ds\right) \right).
\end{equation}
Furthermore, Smarandache curves for $\gamma$ are
\begin{equation*}
\gamma_{\mathbf{e}_{1}\mathbf{e}_{2}}=\left(
\begin{array}{c}
1,-a\int \sinh \left( \int \tau (s)\,ds\right) ds-\sinh \left( \int \tau
(s)\,ds\right) , \\
\cosh \left( \int \tau (s)\,ds\right) +a\int \cosh \left( \int \tau
(s)\,ds\right) \,ds%
\end{array}%
\right),
\end{equation*}
\begin{equation*}
\gamma_{\mathbf{e}_{1}\mathbf{e}_{3}}=\left(
\begin{array}{c}
1,-\cosh \left( \int \tau (s)\,ds\right) -a\int \sinh \left( \int \tau
(s)\,ds\right) ds, \\
a\int \cosh \left( \int \tau (s)\,ds\right) \,ds+\sinh \left( \int \tau
(s)\,ds\right)
\end{array}%
\right),
\end{equation*}
\begin{equation}
\gamma_{\mathbf{e}_{1}\mathbf{e}_{2}\mathbf{e}_{3}}=\left( 1,-e^{\int \tau (s)\,ds}-a\int \sinh \left( \int \tau (s)\,ds\right)
\,ds,e^{\int \tau (s)\,\,ds}+a\int \cosh \left( \int \tau (s)\,ds\right)
\,ds\right).
\end{equation}
\textbf{Case 4.3.2} \ $\gamma(s)$ is \emph{timelike}:
\begin{equation}
\gamma (s)=\left( s,a\int \left( \int \cosh \left( \int \tau (s)\,ds\right)
\,ds\,\right) ds,a\int \left( \int \sinh \left( \int \tau (s)\,ds\right)
\,ds\,\right) ds\right).
\end{equation}
We differentiate this equation three times to get
\begin{equation*}
\gamma ^{\prime }(s)=\left( 1,a\int \cosh \left( \int \tau (s)\,ds\right) \,ds,a\int \sinh \left(
\int \tau (s)\,ds\right) \,ds\right),
\end{equation*}
\begin{equation*}
\gamma ^{\prime \prime }(s)=\left( 0,a\cosh \left( \int \tau (s)\,ds\right) ,a\sinh \left( \int \tau
(s)\,ds\right) \right),
\end{equation*}
\begin{equation}
\gamma ^{\prime \prime \prime }(s)=\left( 0,a\sinh \left( \int \tau (s)\,ds\right) \tau (s),a\cosh \left( \int
\tau (s)\,ds\right) \tau (s)\right).
\end{equation}
The tangent, principal normal and binormal vectors of $\gamma$ are in the following forms
\begin{equation*}
(\mathbf{e}_{1})_{\gamma}=\left( 1,a\int \cosh \left( \int \tau (s)\,ds\right) \,ds,a\int \sinh \left(
\int \tau (s)\,ds\right) \,ds\right),
\end{equation*}
\begin{equation*}
(\mathbf{e}_{2})_{\gamma}=\left( 0,\cosh \left( \int \tau (s)\,ds\right) ,\sinh \left( \int \tau
(s)\,ds\right) \right),
\end{equation*}
\begin{equation}
(\mathbf{e}_{3})_{\gamma}=\left( 0,\sinh \left( \int \tau (s)\,ds\right) ,\cosh \left( \int \tau
(s)\,ds\right) \right).
\end{equation}
So, Smarandache curves for $\gamma$ are as follows
\begin{equation*}
\gamma_{\mathbf{e}_{1}\mathbf{e}_{2}}=\left(
\begin{array}{c}
1,a\int \cosh \left( \int \tau (s)\,ds\right) ds+\cosh \left( \int \tau
(s)\,ds\right) , \\
\sinh \left( \int \tau (s)\,ds\right) +a\int \sinh \left( \int \tau
(s)\,ds\right) \,ds%
\end{array}%
\right),
\end{equation*}
\begin{equation*}
\gamma_{\mathbf{e}_{1}\mathbf{e}_{3}}=\left(
\begin{array}{c}
1,\sinh \left( \int \tau (s)\,ds\right) +a\int \cosh \left( \int \tau
(s)\,ds\right) ds, \\
a\int \sinh \left( \int \tau (s)\,ds\right) \,ds+\cosh \left( \int \tau
(s)\,ds\right)
\end{array}%
\right),
\end{equation*}
\begin{equation}
\gamma_{\mathbf{e}_{1}\mathbf{e}_{2}\mathbf{e}_{3}}=\left( 1,e^{\int \tau (s)\,ds}+a\int \cosh \left( \int \tau (s)\,ds\right)
\,ds,e^{\int \tau (s)\,\,ds}+a\int \sinh \left( \int \tau (s)\,ds\right)
\,ds\right).
\end{equation}

\subsection{Smarandache curves of Anti-Salkowski curve}

Let $\delta (s)$ be Anti-Salkowski curve in $G_{3}^{1}$ with ($\kappa=\kappa (s), \tau=a=const.$)\\
\textbf{Case 4.4.1} \ $\delta(s)$ is \emph{spacelike}:
\begin{equation}
\delta (s)=\left( s,-\int \left( \int \sinh (bs)\kappa (s)\,ds\right) \,ds,\int \left(
\int \cosh (bs)\kappa (s)\,ds\right) \,ds\right).
\end{equation}
It gives us the following derivatives
\begin{equation*}
\delta ^{\prime }(s)=\left( 1,-\int \sinh (bs)\kappa (s)\,\,ds,\int \cosh (bs)\kappa (s)ds\right),
\end{equation*}
\begin{equation*}
\delta ^{\prime \prime }(s)=\left( 0,-\sinh (bs)\kappa (s),\cosh (bs)\kappa (s)\right),
\end{equation*}
\begin{equation}
\delta ^{\prime \prime \prime }(s)=\left( 0,-\kappa ^{\prime }(s)\sinh (bs)-b\cosh (bs)\kappa (s),\kappa
^{\prime }(s)\cosh (bs)+b\sinh (bs)\kappa (s)\right).
\end{equation}
Further, we obtain the following Frenet vectors $\mathbf{e}_{1}$, $\mathbf{e}_{2}$, $\mathbf{e}_{3}$ in the form
\begin{equation*}
(\mathbf{e}_{1})_{\delta}=\left( 1,-\int \sinh (bs)\kappa (s)\,ds,\int \cosh (bs)\kappa (s)\,ds\right),
\end{equation*}
\begin{equation*}
(\mathbf{e}_{2})_{\delta}=\left( 0,-\sinh (bs),\cosh (bs)\right),
\end{equation*}
\begin{equation}
(\mathbf{e}_{3})_{\delta}=\left( 0,-\cosh (bs),\sinh (bs)\right).
\end{equation}
Thus, the above computations of Frenet vectors give Smarandache curves as follows
\begin{equation*}
\delta_{\mathbf{e}_{1}\mathbf{e}_{2}}=\left( 1,-\int \sinh (bs)\kappa (s)\,ds-\sinh (bs),\cosh (bs)+\int \cosh
(bs)\kappa (s)ds\right),
\end{equation*}
\begin{equation*}
\delta_{\mathbf{e}_{1}\mathbf{e}_{3}}=\left( 1,-\cosh (bs)-\int \sinh (bs)\kappa (s)ds,\int \cosh (bs)\kappa
(s)\,ds+\sinh (bs)\right),
\end{equation*}
\begin{equation}
\delta_{\mathbf{e}_{1}\mathbf{e}_{2}\mathbf{e}_{3}}=\left( 1,-e^{bs}-\int \sinh (bs)\kappa (s)ds,e^{bs}+\int \cosh (bs)\kappa
(s)\,ds\right).
\end{equation}\\
\textbf{Case 4.4.2} \ $\delta(s)$ is \emph{timelike}:
\begin{equation}
\delta (s)=\left( s,\int \left( \int \cosh (bs)\kappa (s)\,ds\right) \,ds,\int \left(
\int \sinh (bs)\kappa (s)\,ds\right) \,ds\right).
\end{equation}
The derivatives of $\delta$ are
\begin{equation*}
\delta ^{\prime }(s)=\left( 1,\int \cosh (bs)\kappa (s)\,\,ds,\int \sinh (bs)\kappa (s)ds\right),
\end{equation*}
\begin{equation*}
\delta ^{\prime \prime }(s)=\left( 0,\cosh (bs)\kappa (s),\sinh (bs)\kappa (s)\right),
\end{equation*}
\begin{equation}
\delta ^{\prime \prime \prime }(s)=\left( 0,\kappa ^{\prime }(s)\cosh (bs)+b\sinh (bs)\kappa (s),\kappa
^{\prime }(s)\sinh (bs)+b\cosh (bs)\kappa (s)\right).
\end{equation}
Hence, we obtain the following Frenet vectors
\begin{equation*}
(\mathbf{e}_{1})_{\delta}=\left( 1,\int \cosh (bs)\kappa (s)\,ds,\int \sinh (bs)\kappa (s)\,ds\right),
\end{equation*}
\begin{equation*}
(\mathbf{e}_{2})_{\delta}=\left( 0,\cosh (bs),\sinh (bs)\right),
\end{equation*}
\begin{equation}
(\mathbf{e}_{3})_{\delta}=\left( 0,\sinh (bs),\cosh (bs)\right).
\end{equation}
Thus the Smarandache curves by
\begin{equation*}
\delta_{\mathbf{e}_{1}\mathbf{e}_{2}}=\left( 1,\cosh (bs)+\int \cosh (bs)\kappa (s)ds,\int \sinh (bs)\kappa
(s)\,\,ds+\sinh (bs)\right),
\end{equation*}
\begin{equation*}
\delta_{\mathbf{e}_{1}\mathbf{e}_{3}}=\left( 1,\int \cosh (bs)\kappa (s)ds+\sinh (bs),\cosh (bs)+\int \sinh
(bs)\kappa (s)\,ds\right),
\end{equation*}
\begin{equation}
\delta_{\mathbf{e}_{1}\mathbf{e}_{2}\mathbf{e}_{3}}=\left( 1,e^{bs}+\int \cosh (bs)\kappa (s)\,ds,e^{bs}+\int \sinh (bs)\kappa
(s)ds\right).
\end{equation}
\begin{notation}
In the light of the above calculations, there are not $\mathbf{e}_{2}\mathbf{e}_{3}$ Smarandache curves in the Galilean or pseudo-Galilean spaces.
\end{notation}
\section{Examples}
\begin{example}
Consider $\alpha(u)$ is a \emph{spacelike} general helix in $G_{3}^{1}$ parameterized by
\begin{equation*}
\alpha (s)=\left( u,\frac{1}{6}u\left( -\cosh \left( 2\ln (u)\right) +2\sinh \left(
2\ln (u)\right) \right) ,\frac{1}{6}u\left( 2\cosh \left( 2\ln (u)\right)
-\sinh \left( 2\ln (u)\right) \right) \right).
\end{equation*}
We use the derivatives of $\alpha ; \alpha ^{\prime }, \alpha ^{\prime \prime }, \alpha ^{\prime \prime \prime }$ to get the associated trihedron of $\alpha$ as follows
\begin{equation*}
(\mathbf{e}_{1})_{\alpha}=\left( 1,\frac{1}{2}\cosh \left( 2\ln (u)\right) ,\frac{1}{2}\sinh \left(
2\ln (u)\right) \right),
\end{equation*}
\begin{equation*}
(\mathbf{e}_{2})_{\alpha}=\left( 0,\sinh \left( 2\ln (u)\right) ,\cosh \left( 2\ln (u)\right) \right),
\end{equation*}
\begin{equation*}
(\mathbf{e}_{3})_{\alpha}=\left( 0,-\cosh \left( 2\ln (u)\right) ,-\sinh \left( 2\ln (u)\right)
\right).
\end{equation*}
Curvature $\kappa(s)$ and torsion $\tau(s)$ are obtained as follows
\begin{equation*}
\kappa_{\alpha}=\frac{1}{u}, \tau_{\alpha}=\frac{-2}{u}.
\end{equation*}
According to the above calculations, Smarandache curves of $\alpha$ are
\begin{equation*}
\alpha_{\mathbf{e}_{1}\mathbf{e}_{2}}=\left( 1,\frac{1}{2}\cosh \left( 2\ln (u)\right) +\sinh \left( 2\ln
(u)\right) ,\frac{1+3u^{4}}{4u^{2}}\right),
\end{equation*}
\begin{equation*}
\alpha_{\mathbf{e}_{1}\mathbf{e}_{3}}=\left( 1,-\frac{1}{2}\cosh \left( 2\ln (u)\right) ,-\frac{1}{2}\sinh \left(
2\ln (u)\right) \right),
\end{equation*}
\begin{equation*}
\alpha_{\mathbf{e}_{1}\mathbf{e}_{2}\mathbf{e}_{3}}=\left( 1,\frac{-3+u^{4}}{4u^{2}},\frac{3+u^{4}}{4u^{2}}\right).
\end{equation*}
\end{example}
\begin{center}
\begin{figure}[h]
\centering
\includegraphics[scale=0.7]{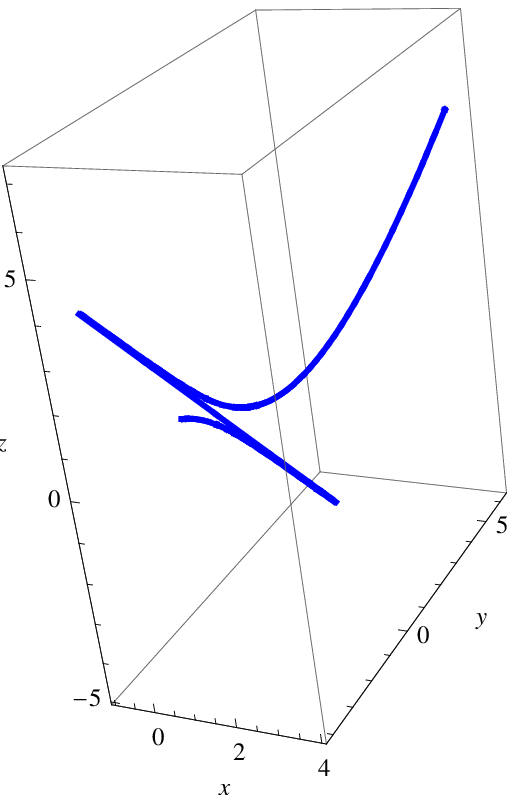}
\label{fig:as}
\caption{The spacelike general helix $\alpha$ in $G_{3}^{1}$ with $\kappa_{\alpha}=\frac{1}{u}$ and $\tau_{\alpha}=\frac{-2}{u}$.}
\end{figure}
\end{center}
\begin{center}
\begin{figure}[h]
\begin{minipage}[t]{0.3\textwidth}
\hspace{-0.1\textwidth}
\centering
\includegraphics[scale=0.6]{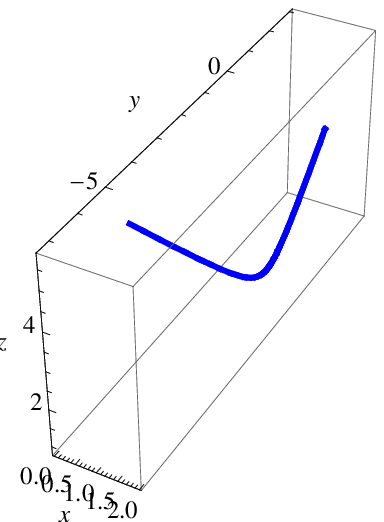}
\label{fig:tnas}
\end{minipage}
\begin{minipage}[t]{0.3\textwidth}
\hspace{-0.1\textwidth}
\centering
\includegraphics[scale=0.6]{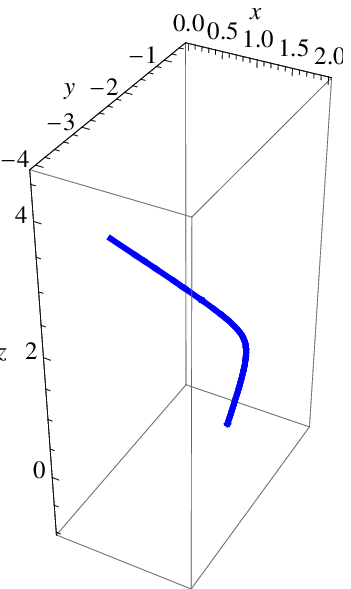}
\label{fig:tbas}
\end{minipage}
\begin{minipage}[t]{0.3\textwidth}
\hspace{-0.1\textwidth}
\centering
\includegraphics[scale=0.6]{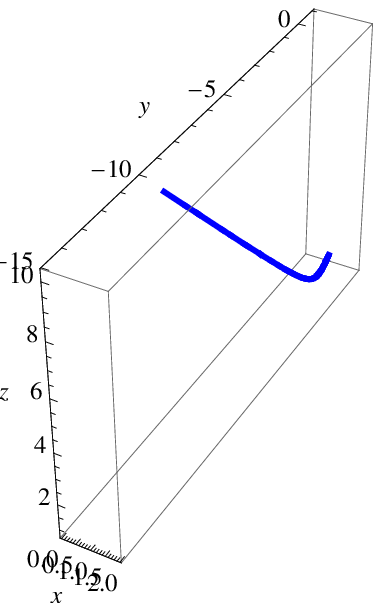}
\label{fig:tnbas}
\end{minipage}
\caption{The $\textbf{e}_{1}\textbf{e}_{2}$, $\textbf{e}_{1}\textbf{e}_{3}$ and $\textbf{e}_{1}\textbf{e}_{2}\textbf{e}_{3}$ Smarandache curves of $\alpha$.}
\end{figure}
\end{center}
\begin{example}
Consider $\alpha^{*} (s)$ is a \emph{timelike} general helix in $G_{3}^{1}$ given by
\begin{equation*}
\alpha^{*} (s)=\left( u,\frac{1}{6}u\left( 2\cosh \left( 2\ln (u)\right) -\sinh \left( 2\ln
(u)\right) \right) ,\frac{1}{6}u\left( -\cosh \left( 2\ln (u)\right) +2\sinh
\left( 2\ln (u)\right) \right) \right).
\end{equation*}
Also, we use the derivatives of $\alpha^{*}; (\alpha^{*}) ^{\prime }, (\alpha^{*}) ^{\prime \prime }, (\alpha^{*}) ^{\prime \prime \prime }$ to get the associated trihedron of $\alpha^{*}$ as follows
\begin{equation*}
(\mathbf{e}_{1})_{\alpha^{*}}=\left( 1,\frac{1}{2}\sinh \left( 2\ln (u)\right) ,\frac{1}{2}\cosh \left(
2\ln (u)\right) \right),
\end{equation*}
\begin{equation*}
(\mathbf{e}_{2})_{\alpha^{*}}=\left( 0,\cosh \left( 2\ln (u)\right) ,\sinh \left( 2\ln (u)\right) \right),
\end{equation*}
\begin{equation*}
(\mathbf{e}_{3})_{\alpha^{*}}=\left( 0,\sinh \left( 2\ln (u)\right) ,\cosh \left( 2\ln (u)\right) \right).
\end{equation*}
Curvature $\kappa(s)$ and torsion $\tau(s)$ are obtained as follows
\begin{equation*}
\kappa_{\alpha^{*}}=\frac{1}{u}, \tau_{\alpha^{*}}=\frac{2}{u}.
\end{equation*}
According to the above calculations, Smarandache curves of $\alpha^{*}$ are
\begin{equation*}
\alpha^{*}_{\mathbf{e}_{1}\mathbf{e}_{2}}=\left( 1,\frac{1+3u^{4}}{4u^{2}},\frac{1}{2}\cosh \left( 2\ln (u)\right)
+\sinh \left( 2\ln (u)\right) \right),
\end{equation*}
\begin{equation*}
\alpha^{*}_{\mathbf{e}_{1}\mathbf{e}_{3}}=\left( 1,\frac{3}{2}\sinh \left( 2\ln (u)\right) ,\frac{3}{2}\cosh \left(
2\ln (u)\right) \right),
\end{equation*}
\begin{equation*}
\alpha^{*}_{\mathbf{e}_{1}\mathbf{e}_{2}\mathbf{e}_{3}}=\left( 1,\cosh \left( 2\ln (u)\right) +\frac{3}{2}\sinh \left( 2\ln
(u)\right) ,\frac{1+5u^{4}}{4u^{2}}\right).
\end{equation*}
\end{example}
\begin{center}
\begin{figure}[h]
\centering
\includegraphics[scale=0.7]{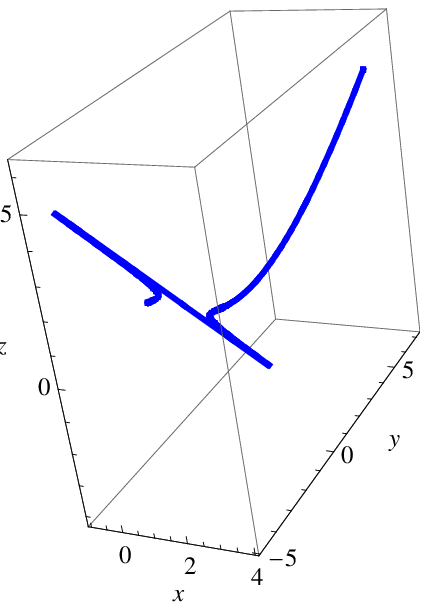}
\label{fig:as}
\caption{The timelike general helix $\alpha^{*}$ in $G_{3}^{1}$ with $\kappa_{\alpha^{*}}=\frac{1}{u}$ and $\tau_{\alpha^{*}}=\frac{2}{u}$.}
\end{figure}
\end{center}
\begin{center}
\begin{figure}[h]
\begin{minipage}[t]{0.3\textwidth}
\hspace{-0.1\textwidth}
\centering
\includegraphics[scale=0.6]{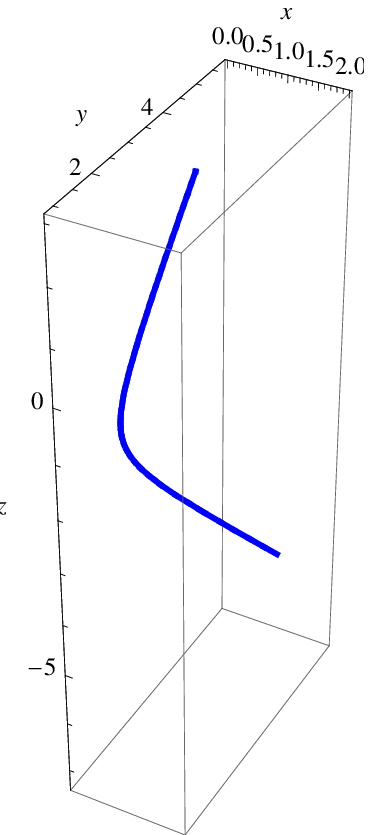}
\label{fig:tnas}
\end{minipage}
\begin{minipage}[t]{0.3\textwidth}
\hspace{-0.1\textwidth}
\centering
\includegraphics[scale=0.6]{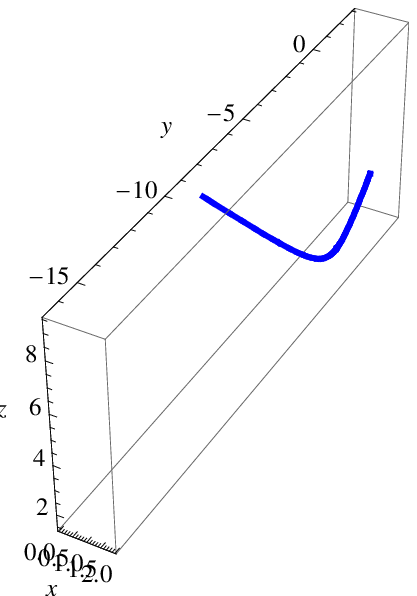}
\label{fig:tbas}
\end{minipage}
\begin{minipage}[t]{0.3\textwidth}
\hspace{-0.1\textwidth}
\centering
\includegraphics[scale=0.6]{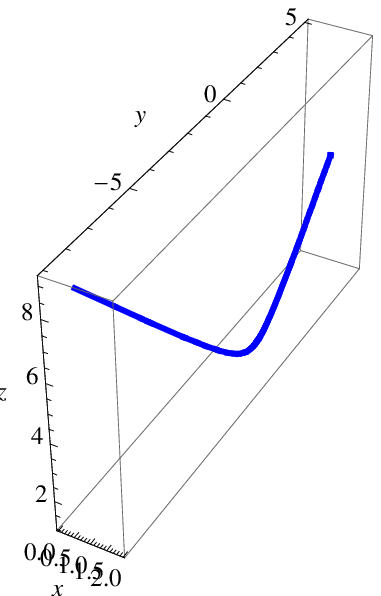}
\label{fig:tnbas}
\end{minipage}
\caption{The $\textbf{e}_{1}\textbf{e}_{2}$, $\textbf{e}_{1}\textbf{e}_{3}$ and $\textbf{e}_{1}\textbf{e}_{2}\textbf{e}_{3}$ Smarandache curves of $\alpha^{*}$.}
\end{figure}
\end{center}
\newpage
\begin{example}
Let $\delta :I\longrightarrow G_{3}^{1}$ be a \emph{spacelike} Anti-Salkowski curve parameterized by
\begin{equation*}
\delta (s)=\left( u,\frac{1}{9}e^{-u}\left( 4\cosh (2u)+5\sinh (2u)\right) ,\frac{1}{9}%
e^{-u}\left( 5\cosh (2u)+4\sinh (2u)\right) \right).
\end{equation*}
By differentiation, we get
\begin{equation*}
\delta ^{\prime }(s)=\left( 1,\frac{1}{6}\left( e^{-3u}+3e^{u}\right) ,-\frac{1}{6}e^{-3u}+\frac{%
e^{u}}{2}\right),
\end{equation*}
\begin{equation*}
\delta ^{\prime \prime }(s)=\left( 0,e^{-u}\sinh (2u),e^{-u}\cosh (2u)\right),
\end{equation*}
\begin{equation*}
\delta ^{\prime \prime \prime }(s)=\left( 0,\frac{1}{2}\left( 3e^{-3u}+e^{u}\right) ,\frac{1}{2}\left(
-3e^{-3u}+e^{u}\right) \right).
\end{equation*}
Using $(2.5)$ to obtain
\begin{equation*}
(\mathbf{e}_{1})_{\delta}=\left( 1,\frac{1}{6}\left( e^{-3u}+3e^{u}\right) ,-\frac{1}{6}e^{-3u}+\frac{%
e^{u}}{2}\right),
\end{equation*}
\begin{equation*}
(\mathbf{e}_{2})_{\delta}=\left( 0,\sinh (2u),\cosh (2u)\right),
\end{equation*}
\begin{equation*}
(\mathbf{e}_{3})_{\delta}=\left( 0,-\cosh (2u),-\sinh (2u)\right).
\end{equation*}
The natural equations of this curve are given by
\begin{equation*}
\kappa_{\delta}=e^{-u}, \tau_{\delta}=-2.
\end{equation*}
Thus, the Smarandache curves of $\delta$ are respectively, given by
\begin{equation*}
\delta_{\mathbf{e}_{1}\mathbf{e}_{2}}=\left( 1,\frac{1}{6}\left( e^{-3u}+3e^{u}\right) +\sinh (2u),-\frac{1}{6}%
e^{-3u}+\frac{e^{u}}{2}+\cosh (2u)\right),
\end{equation*}
\begin{equation*}
\delta_{\mathbf{e}_{1}\mathbf{e}_{3}}=\left( 1,\frac{1}{6}\left( e^{-3u}+3e^{u}-6\cosh (2u)\right) ,-\frac{1}{6}%
e^{-3u}+\frac{e^{u}}{2}-\sinh (2u)\right),
\end{equation*}
\begin{equation*}
\delta_{\mathbf{e}_{1}\mathbf{e}_{2}\mathbf{e}_{3}}=\left( 1,\frac{1}{6}e^{-3u}\left( 1-6e^{u}+3e^{4u}\right) ,\frac{1}{6}%
e^{-3u}\left( -1+6e^{u}+3e^{4u}\right) \right).
\end{equation*}
\end{example}
\begin{center}
\begin{figure}[h]
\centering
\includegraphics[scale=0.7]{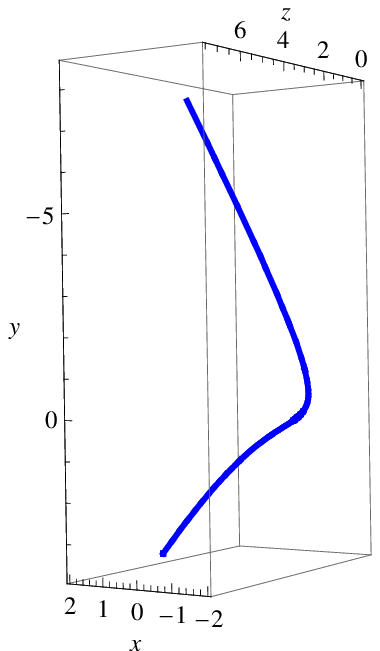}
\label{fig:gh}
\caption{The spacelike Anti-Salkowski curve $\delta$ in $G_{3}^{1}$ with $\kappa_{\delta} =e^{-u}$ and $\tau_{\delta} =-2$.}
\end{figure}
\end{center}
\begin{center}
\begin{figure}[h]
\begin{minipage}[t]{0.3\textwidth}
\hspace{-0.1\textwidth}
\centering
\includegraphics[scale=0.7]{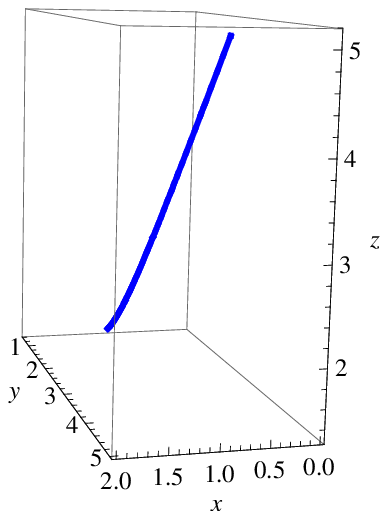}
\label{fig:tngh}
\end{minipage}
\begin{minipage}[t]{0.3\textwidth}
\hspace{-0.1\textwidth}
\centering
\includegraphics[scale=0.7]{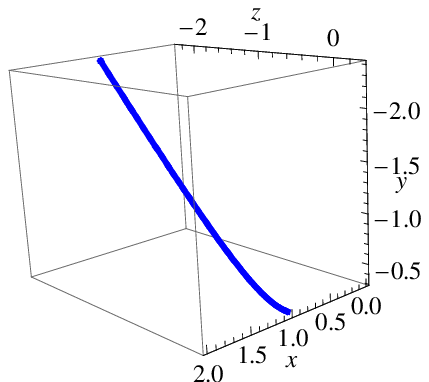}
\label{fig:tbgh}
\end{minipage}
\begin{minipage}[t]{0.3\textwidth}
\hspace{-0.1\textwidth}
\centering
\includegraphics[scale=0.7]{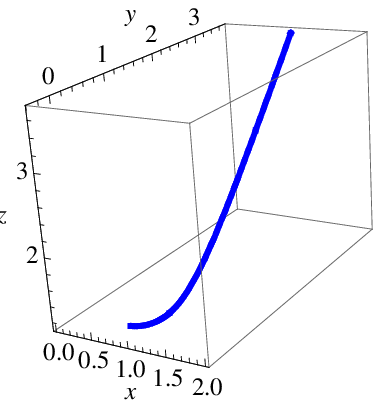}
\label{fig:tnbgh}
\end{minipage}
\caption{From left to right, the $\textbf{e}_{1}\textbf{e}_{2}$, $\textbf{e}_{1}\textbf{e}_{3}$ and $\textbf{e}_{1}\textbf{e}_{2}\textbf{e}_{3}$ Smarandache curves of $\delta$.}
\end{figure}
\end{center}
\begin{example}
Let $\delta^{*}$ be a \emph{timelike} Anti-Salkowski curve in $G_{3}^{1}$ given by
\begin{equation*}
\delta^{*} (s)=\left( u,\frac{1}{9}e^{-u}\left( 5\cosh (2u)+4\sinh (2u)\right) ,\frac{1}{9}%
e^{-u}\left( 4\cosh (2u)+5\sinh (2u)\right) \right).
\end{equation*}
By differentiation, we get
\begin{equation*}
(\delta^{*}) ^{\prime }(s)=\left( 1,-\frac{1}{6}e^{-3u}+\frac{e^{u}}{2},\frac{1}{6}\left(
e^{-3u}+3e^{u}\right) \right),
\end{equation*}
\begin{equation*}
(\delta^{*}) ^{\prime \prime }(s)=\left( 0,e^{-u}\cosh (2u),e^{-u}\sinh (2u)\right),
\end{equation*}
\begin{equation*}
(\delta^{*}) ^{\prime \prime \prime }(s)=\left( 0,\frac{1}{2}\left( -3e^{-3u}+e^{u}\right) ,\frac{1}{2}\left(
3e^{-3u}+e^{u}\right) \right).
\end{equation*}
Using $(2.5)$ to obtain
\begin{equation*}
(\mathbf{e}_{1})_{\delta^{*}}=\left( 1,-\frac{1}{6}e^{-3u}+\frac{e^{u}}{2},\frac{1}{6}\left(
e^{-3u}+3e^{u}\right) \right),
\end{equation*}
\begin{equation*}
(\mathbf{e}_{2})_{\delta^{*}}=\left( 0,\cosh (2u),\sinh (2u)\right),
\end{equation*}
\begin{equation*}
(\mathbf{e}_{3})_{\delta^{*}}=\left( 0,\sinh (2u),\cosh (2u)\right).
\end{equation*}
The natural equations of this curve are given by
\begin{equation*}
\kappa_{\delta^{*}}=e^{-u}, \tau_{\delta^{*}}=2.
\end{equation*}
Thus, the Smarandache curves of $\delta^{*}$ are respectively, given by
\begin{equation*}
(\delta^{*})_{\mathbf{e}_{1}\mathbf{e}_{2}}=\left( 1,-\frac{1}{6}e^{-3u}+\frac{e^{u}}{2}+\cosh (2u),\frac{1}{6}\left(
e^{-3u}+3e^{u}\right) +\sinh (2u)\right),
\end{equation*}
\begin{equation*}
(\delta^{*})_{\mathbf{e}_{1}\mathbf{e}_{3}}=\left( 1,-\frac{1}{6}e^{-3u}+\frac{e^{u}}{2}+\sinh (2u),\frac{1}{6}\left(
e^{-3u}+3e^{u}\right) +\cosh (2u)\right),
\end{equation*}
\begin{equation*}
(\delta^{*})_{\mathbf{e}_{1}\mathbf{e}_{2}\mathbf{e}_{3}}=\left( 1,-\frac{1}{6}e^{-3u}+\frac{e^{u}}{2}+e^{2u},\frac{e^{-3u}}{6}+\frac{
e^{u}}{2}+e^{2u}\right).
\end{equation*}
\end{example}
\begin{center}
\begin{figure}[h]
\centering
\includegraphics[scale=0.7]{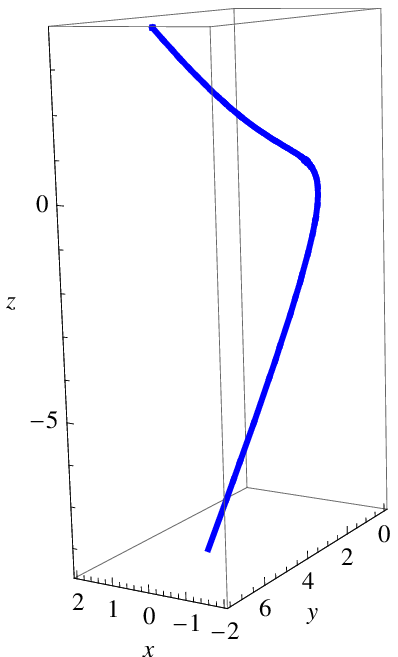}
\label{fig:gh}
\caption{The timelike Anti-Salkowski curve $\delta^{*}$ in $G_{3}^{1}$ with $\kappa_{\delta^{*}} =e^{-u}$ and $\tau_{\delta^{*}} =2$.}
\end{figure}
\end{center}
\begin{center}
\begin{figure}[h]
\begin{minipage}[t]{0.3\textwidth}
\hspace{-0.1\textwidth}
\centering
\includegraphics[scale=0.7]{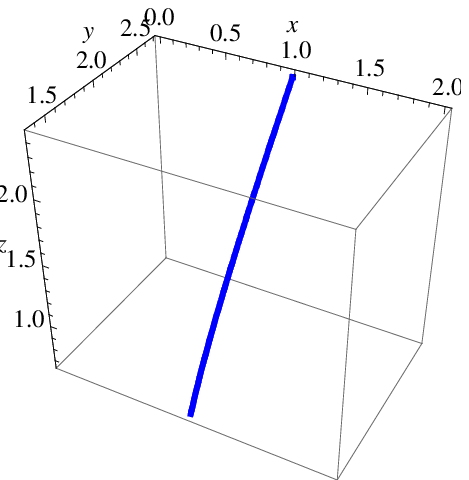}
\label{fig:tngh}
\end{minipage}
\begin{minipage}[t]{0.3\textwidth}
\hspace{-0.1\textwidth}
\centering
\includegraphics[scale=0.7]{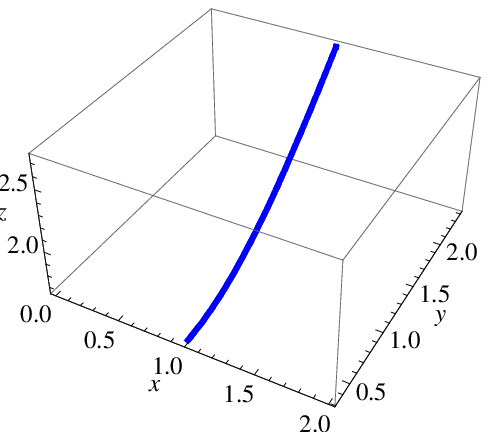}
\label{fig:tbgh}
\end{minipage}
\begin{minipage}[t]{0.3\textwidth}
\hspace{-0.1\textwidth}
\centering
\includegraphics[scale=0.7]{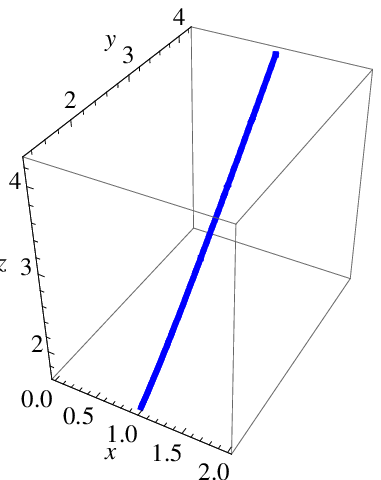}
\label{fig:tnbgh}
\end{minipage}
\caption{From left to right, the $\textbf{e}_{1}\textbf{e}_{2}$, $\textbf{e}_{1}\textbf{e}_{3}$ and $\textbf{e}_{1}\textbf{e}_{2}\textbf{e}_{3}$ Smarandache curves of $\delta^{*}$.}
\end{figure}
\end{center}
\begin{example}
Consider $\eta$ is a \emph{timelike} spiral in $G_{3}^{1}$ parameterized as follows
\begin{equation*}
\eta (s)=\left( u,(2+u)(-1+\ln (2+u)),0\right).
\end{equation*}
So, we get
\begin{equation*}
\eta ^{\prime }(s)=\left( 1,\ln (2+u),0\right),
\end{equation*}
\begin{equation*}
\eta ^{\prime \prime }(s)=\left( 0,\frac{1}{2+u},0\right),
\end{equation*}
\begin{equation*}
\eta ^{\prime \prime \prime }(s)=\left( 0,-\frac{1}{(2+u)^{2}},0\right).
\end{equation*}
The Frenet vectors of $\eta$ are
\begin{equation*}
(\mathbf{e}_{1})_{\eta}=\left( 1,\ln (2+u),0\right),
\end{equation*}
\begin{equation*}
(\mathbf{e}_{2})_{\eta}=\left( 0,1,0\right),
\end{equation*}
\begin{equation*}
(\mathbf{e}_{3})_{\eta}=\left( 0,0,1\right).
\end{equation*}
The curvatures of this curve are given by
\begin{equation*}
\kappa_{\eta}=\frac{1}{2+u}, \tau_{\eta}=0.
\end{equation*}
Thus, the Smarandache curves of this spiral are given by
\begin{equation*}
\eta_{\mathbf{e}_{1}\mathbf{e}_{2}}=\left( 1,\ln (2+u),1\right),
\end{equation*}
\begin{equation*}
\eta_{\mathbf{e}_{1}\mathbf{e}_{3}}=\left( 1,1+\ln (2+u),1\right),
\end{equation*}
\begin{equation*}
\eta_{\mathbf{e}_{1}\mathbf{e}_{2}\mathbf{e}_{3}}=\left( 1,1+\ln (2+u),1\right).
\end{equation*}
\end{example}
\begin{center}
\begin{figure}[h]
\centering
\includegraphics[scale=0.7]{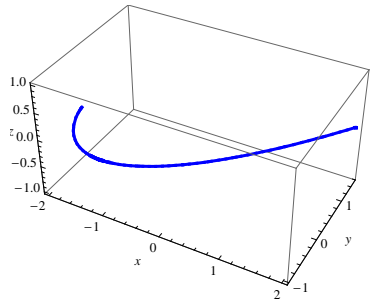}
\label{fig:gh}
\caption{The timelike spiral curve $\eta$ in $G_{3}^{1}$ with $\kappa_{\eta} =\frac{1}{2+u}$ and $\tau_{\eta} =0$.}
\end{figure}
\end{center}
\begin{center}
\begin{figure}[h]
\begin{minipage}[t]{0.3\textwidth}
\hspace{-0.1\textwidth}
\centering
\includegraphics[scale=0.7]{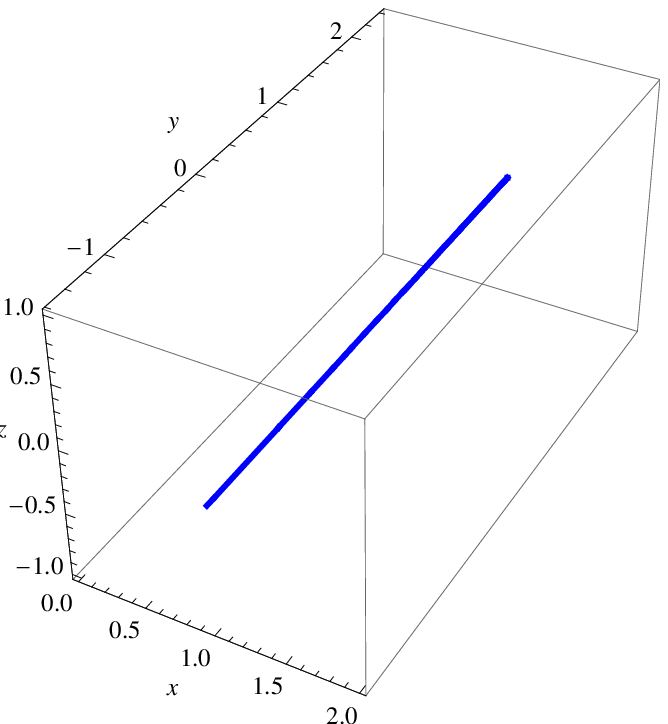}
\label{fig:tngh}
\end{minipage}
\begin{minipage}[t]{0.3\textwidth}
\hspace{-0.1\textwidth}
\centering
\includegraphics[scale=0.7]{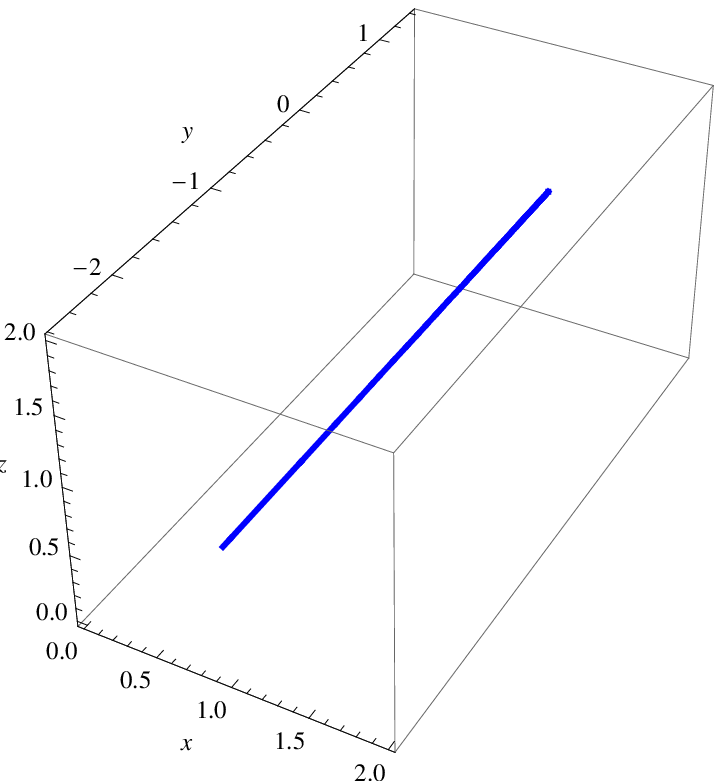}
\label{fig:tbgh}
\end{minipage}
\begin{minipage}[t]{0.3\textwidth}
\hspace{-0.1\textwidth}
\centering
\includegraphics[scale=0.7]{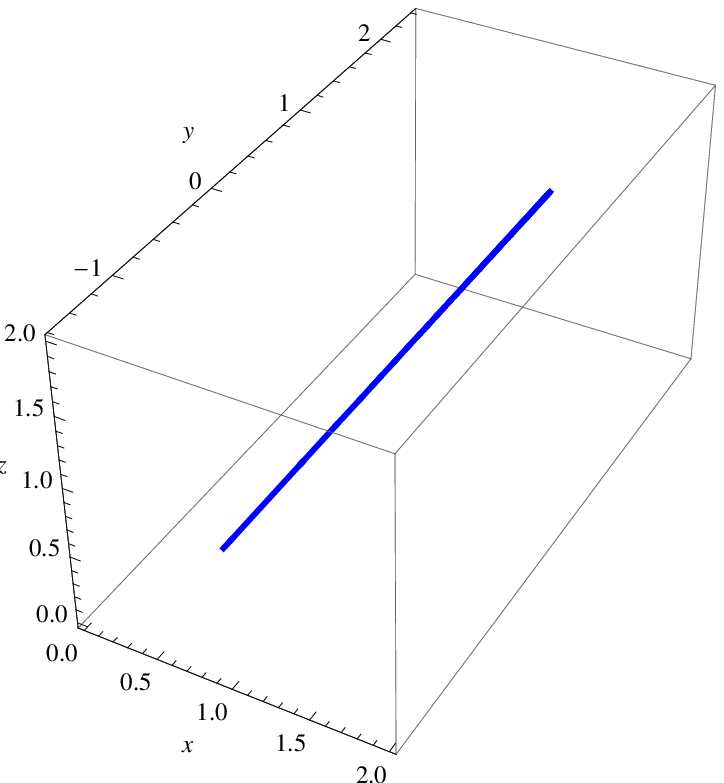}
\label{fig:tnbgh}
\end{minipage}
\caption{The $\textbf{e}_{1}\textbf{e}_{2}$, $\textbf{e}_{1}\textbf{e}_{3}$ and $\textbf{e}_{1}\textbf{e}_{2}\textbf{e}_{3}$ Smarandache curves of $\eta$.}
\end{figure}
\end{center}

\section{Conclusion}

In the three-dimensional pseudo-Galilean space $G_{3}^{1}$, Smarandache curves of space and timelike arbitrary curve and some of its special curves have been studied. Some examples of these curves such as  general helix, Ant-Salkowski and spiral curves have been given and plotted.

\end{document}